\documentclass{amsart}

\newtheorem{theorem}{Theorem}
\newtheorem{corollary}[theorem]{Corollary}
\newtheorem{lemma}[theorem]{Lemma}

\newcommand{\RR}{\mathbb{R}}
\newcommand{\NN}{\mathbb{N}}
\newcommand{\calA}{\mathcal{A}}
\newcommand{\calB}{\mathcal{B}}
\newcommand{\calT}{\mathcal{T}}

\newcommand{\ra}{\rightarrow}

\newcommand{\com}{\mathbin{{\scriptstyle \circ }}}

\newcommand{\sign}{\mathord{\mathrm{sign}}}
\newcommand{\supp}{\mathord{\mathrm{supp}}}
\newcommand{\germ}{\mathord{\mathrm{germ}}}
\newcommand{\Crs}[1]{\mathord{\mathit{C}^{r}_{#1}}}
\newcommand{\Cr}{\mathord{\mathit{C}^{r}}}
\newcommand{\Cinf}{\mathord{\mathit{C}^{\infty}}}

\begin{document}

\title[Multiplicative bijections]
      {Multiplicative bijections between algebras of differentiable functions}

\author{J. Mr\v{c}un}
\author{P. \v{S}emrl}

\address{Department of Mathematics, University of Ljubljana,
         Jadranska 19, 1000 Ljubljana, Slovenia}
\email{janez.mrcun@fmf.uni-lj.si}
\email{peter.semrl@fmf.uni-lj.si}

\thanks{This work was supported in part by the Slovenian Ministry of Science}

\subjclass[2000]{Primary 58A05; Secondary 46E25}
\date{}

\begin{abstract}
We show that any multiplicative bijection between the algebras
of differentiable functions, defined
on differentiable manifolds of positive dimension,
is an algebra isomorphism, given by composition
with a unique diffeomorphism.
\end{abstract}

\maketitle

\section{Introduction} \label{section1}

In the theory of classical algebras, the problem
of characterizing automorphisms is of
fundamental importance. 
In the case of the algebra of
smooth functions on a smooth manifold,
every automorphism is a
composition operator. 

It has been known for a long time that the linear structure
of an algebra is often completely determined by the
multiplicative one.
Already in 1940, Eidelheit \cite{Eidelheit} observed
that the multiplicative bijective maps on real operator algebras
are automatically linear.
An interested reader can find a pure ring-theoretic result
on automatic additivity of multiplicative maps in \cite{Martindale}.
However, this result is not relevant for the case of commutative
rings. In the commutative case,
the situation is more complicated.
For example, suppose that
$\tau\!:X\ra Y$ is a homeomorphism of compact
Hausdorff spaces, and let $p$ be a positive continuous function on $X$.
Then the map $\calT\!:C(Y)\ra C(X)$ between the
corresponding algebras of real continuous functions,
given by 
\[
\calT(g)(x) = | g(\tau (x)) |^{p(x)} \, \sign (g(\tau (x)))\;,
\;\;\;\;\;\;\;\;\;\; x\in X\;, \;\; g\in C(Y)\;,
\]
is a bijective multiplicative map, non-linear if $p\neq 1$.
It turns out \cite{Milgram} that in the absence of isolated points,
every semigroup isomorphism of $C(Y)$ onto
$C(X)$ is of this form.
This result shows that the multiplicative semigroup structure
of $C(X)$ completely determines the underlying space $X$. Namely, if
$C(X)$ and $C(Y)$ are isomorphic as multiplicative semigroups, then by the
above result the spaces $X$ and $Y$ are homeomorphic.

The aim of this note is to show that the result is much nicer
in the case of real differentiable functions on differentiable
manifolds. Here, any semigroup 
isomorphism is automatically linear. More precisely:

\begin{theorem} \label{theorem1}
Let $M$ and $N$ be Hausdorff $\Cr$-manifolds of positive dimension,
$1\leq r <\infty$.
Then for any multiplicative bijection $\calB\!:\Cr(N)\ra\Cr(M)$
there exists a unique $\Cr$-diffeomorphism $\phi\!:M\ra N$
such that
\[ \calB(g)(x)=g(\phi(x)) \]
for any $g\in\Cr(N)$ and any $x\in M$.
In particular, the map $\calB$ is an algebra isomorphism.
\end{theorem}

It should be mentioned that the manifolds $M$ and $N$
appearing in Theorem \ref{theorem1} are not assumed to be
second-countable, paracompact or connected.
Even in the presence of the
linearity assumption,
the above result has been proved in full
generality only very recently \cite{Grabowski,Mrcun}.
As pointed out by A. Weinstein in 2003,
all the earlier proofs (e.g.\ \cite{Bkouche})
depend heavily on the second-countability
assumption.

Unfortunately our methods do not work for the case $r=\infty$.
However, it is tempting to believe that the result holds in this case
as well.

\section{Proof of the main theorem} \label{section2}

Throughout this paper, we assume that $M$ and $N$ are Hausdorff
$\Cr$-manifolds, not necessarily second-countable,
paracompact or connected, $r=0,1,\ldots,\infty$, and that
\[ \calB\!:\Cr(N)\ra\Cr(M) \]
is a multiplicative bijection between the algebras of real
$\Cr$-functions on $N$, respectively $M$.

First, we will show that $\calB$
induces a homeomorphism $M\ra N$ in a natural way, by means
of characteristic sequences of functions introduced in \cite{Mrcun}.
Recall that a sequence $(f_{i})$ of $\Cr$-functions on $M$
is called {\em characteristic at} $x\in M$ if
\begin{enumerate}
\item [(i)]  $f_{i}f_{i+1}=f_{i+1}$ for any $i$, and
\item [(ii)] the associated sequence of supports
             $(\supp(f_{i}))$ is a fundamental system of neighbourhoods
             of $x$ in $M$.
\end{enumerate}
In particular, if $(f_{i})$ is a characteristic sequence
of $\Cr$-functions on $M$ at $x$, then
$\bigcap_{i}\supp(f_{i})=\{x\}$, 
$f_{i}$ equals $1$ on $\supp(f_{i+1})$
and $\supp(f_{i})$
is compact for any $i$ large enough \cite[Lemma 3]{Mrcun}.

For any point $x\in M$ we can choose a characteristic sequence
$(f_{i})$ of $\Cr$-functions on $M$ at $x$. It follows from the proof
of \cite[Lemma 4]{Mrcun} (which is stated for isomorphisms
of algebras, but its proof only uses the fact that the
map between algebras is a multiplicative bijection)
that the sequence $(\calB^{-1}(f_{i}))$ of $\Cr$-functions on $N$
is characteristic at a point $\phi(x)\in N$, and that
this point is independent on the choice of the
sequence $(f_{i})$. In particular, we obtain a map
\[ \phi\!:M\ra N \;. \]
By symmetry, the same construction can be applied to $\calB^{-1}$,
and the associated map $N\ra M$ is clearly the inverse
of $\phi$. In particular, the map $\phi$ is a bijection.

We would like to show that for any $g\in\Cr(N)$, the value
of $\calB(g)$ at $x\in M$ depends only on the value of $g$ at $\phi(x)$.
Actually, this is the core of our problem, and we will
need to do several small steps before achieving this goal.
First, we will see that the value $\calB(g)(x)$ depends only on 
the germ of $g$ at $\phi(x)$. We shall denote by
$\germ_{\phi(x)}(g)=g_{\phi(x)}$ the germ of $g$ at $\phi(x)$,
and by $\Crs{N}$ the sheaf of germs of all
$\Cr$-functions on $N$.

\begin{lemma} \label{lemma2}
For any $g,h\in\Cr(N)$ and any $x\in M$ we have:
\begin{enumerate}
\item [(i)]   $g_{\phi(x)}=h_{\phi(x)}$ if and only if
              $\calB(g)_{x}=\calB(h)_{x}$
\item [(ii)]  $g_{\phi(x)}=1$ if and only if
              $\calB(g)_{x}=1$
\item [(iii)] $g_{\phi(x)}=-1$ if and only if
              $\calB(g)_{x}=-1$
\item [(iv)]  $g(\phi(x))=0$ if and only if
              $\calB(g)(x)=0$
\item [(v)]   $g(\phi(x))>0$ if and only if
              $\calB(g)(x)>0$      
\item [(vi)]   $g(\phi(x))<0$ if and only if
              $\calB(g)(x)<0$
\end{enumerate}
In particular $\calB(1)=1$, $\calB(-1)=-1$, $\calB(-g)=-\calB(g)$, and $\calB(0)=0$.
\end{lemma}

\begin{proof}
First note that by symmetry it is sufficient to prove
only one of the implications in all of the equivalences above.

(i)
Suppose that $g_{\phi(x)}=h_{\phi(x)}$. 
Choose a characteristic sequence
$(f_{i})$ of $\Cr$-functions on $M$ at $x$,
and put $g_{i}=\calB^{-1}(f_{i})$ for any $i$.
By the definition of $\phi$, the sequence
$(g_{i})$ is characteristic at $\phi(x)$.
In particular, the sequence
$(\supp(g_{i}))$ is a fundamental system of neighbourhoods
of $\phi(x)$, hence there exists $k$ large enough so that
\[ g g_{k} = h g_{k} \;. \]
Since $\calB$ is multiplicative, this implies
\begin{equation} \label{equation1}
\calB(g) f_{k} = \calB(h) f_{k} \;.
\end{equation}
Because $(f_{i})$ is characteristic at $x$, the function
$f_{k}$ equals $1$ on a neighbourhood
of $x$, so (\ref{equation1}) gives $\calB(g)_{x}=\calB(h)_{x}$.

(ii)
The map $\calB$ is a multiplicative bijection, therefore it
preserves the unit, i.e.\ $\calB(1)=1$.
If $g_{\phi(x)}=1$, it follows from (i) that
$\calB(g)_{x}=\calB(1)_{x}=1$.

(iii)
If $g_{\phi(x)}=-1$, then $g_{\phi(x)}^{2}=1$, so by (ii)
we have $\calB(g)_{x}^{2}=\calB(g^{2})_{x}=1$.
Because (ii) implies that $\calB(g)_{x}\neq 1$,
it follows that $\calB(g)_{x}=-1$.

(iv)
If $g(\phi(x))\neq 0$, there exists
$w\in\Cr(N)$ such that $(gw)_{\phi(x)}=g_{\phi(x)}w_{\phi(x)}=1$.
It follows from (i),(ii) and the multiplicativity of $\calB$ that
$\calB(gw)_{x}=\calB(g)_{x}\calB(w)_{x}=1$, and in particular
$\calB(g)(x)\neq 0$.

(v)
If $g(\phi(x))> 0$, there exists
$w\in\Cr(N)$ such that $w^{2}_{\phi(x)}=g_{\phi(x)}$.
It follows that $\calB(w)^{2}_{x}=\calB(w^{2})_{x}=\calB(g)_{x}$,
so $\calB(g)(x)\geq 0$. On the other hand, we know from (iv) that
$\calB(g)(x)\neq 0$.

(vi)
This follows directly from (iv) and (v).
\end{proof}

\begin{lemma} \label{lemma3}
The map $\phi\!:M\ra N$ is a homeomorphism. Furthermore,
the map $\calB$ induces
a multiplicative homeomorphism of sheaves
$\Crs{N}\ra\Crs{M}$ over $\phi^{-1}$,
again denoted by $\calB$,
which is given by
\[ \calB(g_{\phi(x)})=\calB(g)_{x} \]
for any $x\in M$ and $g\in\Cr(N)$.
\end{lemma}

\begin{proof}
Take any open subset $V$ of $N$ and let
$x\in\phi^{-1}(V)$. Choose a function $g\in\Cr(N)$ with
$\supp(g)\subset V$ such that $g(\phi(x))\neq 0$, and put
\[ U=\{ x'\in M \,|\, \calB(g)(x')\neq 0 \} \;. \]
First observe that $U$ is open in $M$ because $\calB(g)$ is continuous.
It follows from Lemma \ref{lemma2} (iv)
that $U\subset \phi^{-1}(V)$
and that $x\in U$. This shows that $\phi^{-1}(V)$ is open.
We therefore conclude that $\phi$ is continuous.

A symmetrical argument shows that
$\phi^{-1}$ is continuous as well, so $\phi$ is a homeomorphism.
The rest of the statement follows from Lemma \ref{lemma2} (i).
\end{proof}

As a consequence of Lemma \ref{lemma3},
the manifolds $M$ and $N$ have the same dimension.
From now on we will assume that $n=\dim M=\dim N\geq 1$.

It also follows from Lemma \ref{lemma3} that
for any open subset $V$ of $N$ we have the
multiplicative bijection
$\calB_{V}\!:\Cr(V)\ra\Cr(\phi^{-1}(V))$
such that
$\calB_{V}(h)_{x}=\calB(h_{\phi(x)})$
for any $h\in\Cr(V)$ and
for any $x\in\phi^{-1}(V)$.
Observe that to prove Theorem \ref{theorem1} it
is sufficient to find an open cover
$(V_{j})$ of $N$ such that 
the multiplicative bijection
$\calB_{V_{j}}$ is given by composition with
$\phi|_{\phi^{-1}(V_{j})}$, for any $j$.
For instance, it is suitable to choose the
cover so that $V_{j}$ and $\phi^{-1}(V_{j})$
are coordinate charts on $N$, respectively $M$.

Define a map $\calA\!:\Cr(N)\ra\Cr(M)$ by
\[ \calA(g) = {\ln} \com \calB( {\exp} \com g ) \;. \]
Note that $\calA$ is an additive bijection, with inverse
$\calA^{-1}(f)={\ln} \com \calB^{-1}( {\exp} \com f )$.
The properties of $\calB$, stated in Lemma \ref{lemma2} and
Lemma \ref{lemma3},
obviously translate into analogous properties of the map $\calA$.

\begin{corollary} \label{corollary4}
For any $g,h\in\Cr(N)$ and any $x\in M$ we have
$g_{\phi(x)}=h_{\phi(x)}$ if and only if
$\calA(g)_{x}=\calA(h)_{x}$.
Furthermore, the map $\calA$ induces an additive homeomorphism
of sheaves $\calA\!:\Crs{N}\ra\Crs{M}$ over $\phi^{-1}$, and extends
to additive bijections $\calA_{V}\!:\Cr(V)\ra\Cr(\phi^{-1}(V))$
by $\calA_{V}(h)_{x}=\calA(h_{\phi(x)})$,
for all open subsets $V$ of the manifold $N$.
\end{corollary}

Let $f$ be a real $\Cr$-function defined on an open subset $U$ of $M$.
Suppose that $(x_{1},\ldots,x_{n})\!:W\ra\RR^{n}$ are local coordinates 
on an open subset $W$ of $U$,
and let $\alpha=(\alpha_{1},\ldots,\alpha_{n})$
be a multi-index of order $|\alpha|=\alpha_{1}+\cdots+\alpha_{n}\leq r$
($\alpha_{i}\in\NN\cup\{0\}$).
We use the standard notation
\[ D^{\alpha}(f) =
\frac{\partial^{|\alpha|}f}
     {\partial x^{\alpha_{1}}_{1} \cdots \partial x^{\alpha_{n}}_{n}} \]
for the partial $\alpha$-derivative of $f$ of order $|\alpha|$ on $W$.
We shall write $j^{k}_{x}(f)$ for the $k$-jet of $f$
at a point $x\in U$, $k=0,1,\ldots,r$.
It is an equivalence class of real $\Cr$-functions defined
on open neighbourhoods of $x$, with two functions belonging
to the same $k$-jet at $x$ if and only if they
have the same partial derivatives at $x$ of orders
$0,1,\ldots,k$ with respect to
some (or any) local coordinates around $x$ (see \cite{Saunders}).

\begin{lemma} \label{lemma5}
If $g,h\in\Cr(N)$ satisfy
$j^{r}_{\phi(x)}(g)=j^{r}_{\phi(x)}(h)$ for some $x\in M$,
then $j^{r}_{x}(\calA(g))=j^{r}_{x}(\calA(h))$.
\end{lemma}

\begin{proof}
By additivity of $\calA$ we can assume without loss of generality
that $h=0$.
Choose a $\Cr$-diffeomorphism $\psi\!:V\ra\RR^{n}$, defined on
an open neighbourhood $V\subset N$ of $\phi(x)$,
such that $\psi(\phi(x))=0$.
The function $g\com\psi^{-1}\in\Cr(\RR^{n})$ satisfies
$j^{r}_{0}(g\com\psi^{-1})=j^{r}_{0}(0)$
because $j^{r}_{\phi(x)}(g)=j^{r}_{\phi(x)}(0)$.
Hence, by the Whitney extension theorem \cite{Whitney},
there exists a function $h\in\Cr(\RR^{n})$ such that
\[
h|_{[-1,0]^{n}} = 0
\]
and
\[
h|_{[0,1]^{n}} = (g\com\psi^{-1})|_{[0,1]^{n}} \;.
\]
Since $g|_{V}=h\com\psi + (g|_{V} - h\com\psi)$,
it follows that
\begin{equation} \label{equation2}
\calA(g)_{x}=\calA(g_{\phi(x)})=\calA((h\com\psi)_{\phi(x)})
             + \calA(g_{\phi(x)}-(h\com\psi)_{\phi(x)}) \;. 
\end{equation}
By construction of $h$ we have
$h|_{[-1,0]^{n}}=0$ and
$(g\com\psi^{-1}-h)|_{[0,1]^{n}}=0$.
Therefore, Corollary \ref{corollary4} implies that
\[
\calA_{V}(h\com\psi)_{\phi^{-1}(\psi^{-1}(u))}=0
\]
and
\[
\calA_{V}(g|_{V}-h\com\psi)_{\phi^{-1}(\psi^{-1}(v))}=0
\]
for any $u\in (-1,0)^{n}$ and any
$v\in (0,1)^{n}$.
In particular,
all the partial derivatives of order $0,1,\ldots,r$ of
$\calA_{V}(h\com\psi)$ and $\calA_{V}(g|_{V}-h\com\psi)$
are zero arbitrary close to $x$. Since these two
are both $\Cr$-functions on a neighbourhood of $x$,
this implies
\[ j^{r}_{x}(\calA_{V}(h\com\psi))=j^{r}_{x}(0) \]
and
\[ j^{r}_{x}(\calA_{V}(g|_{V}-h\com\psi))=j^{r}_{x}(0) \;, \]
therefore $j^{r}_{x}(\calA(g))=j^{r}_{x}(0)$ by (\ref{equation2}).
\end{proof}

\begin{lemma} \label{lemma6}
Suppose that $r<\infty$ and that
$M$ and $N$ are open subsets of $\RR^{n}$.
There exist an open subset $R$ with
discrete complement in $M$
and continuous real functions $p_{\alpha}$ on $R$,
for any multi-index $\alpha$ of order $|\alpha|\leq r$,
such that
\[ \calA(g)(x) = \sum_{|\alpha|\leq r}
                 p_{\alpha}(x) D^{\alpha}(g)(\phi(x)) \]
for any $g\in \Cr(N)$ and any $x\in R$.
\end{lemma}

\begin{proof}
Denote by $\Pi_{r}=\RR_{r}[t_{1},\ldots,t_{n}]$ the finite dimensional
subspace of $\Cr(N)$
of polynomials of order at most $r$.
Recall that any polynomial in $\Pi_{r}$ is uniquely determined by
the values of all its partial derivatives of
order at most $r$ at any fixed point $y\in N$,
and that this parameterization of $\Pi_{r}$
is linear.

For any $x\in M$ define an additive map
$\Phi_{x}\!:\Pi_{r}\ra \RR$ by
\[ \Phi_{x}(P)=\calA(P)(x)\;, \]
and let
\[ R = \{ x\in M \,|\, \Phi_{x}
     \textrm{ is bounded on a neighbourhood of } 0\in \Pi_{r} \} \;. \]
First, we will show that $M\setminus R$ is closed and discrete in $M$.
To this end, suppose that
$(x_{i})$ is an injective sequence of points in $M\setminus R$
which converges to $x\in M$. 
Put $y_{i}=\phi(x_{i})$, choose a positive decreasing sequence
$(\epsilon_{i})$ converging to $0$ such that the open balls
$K(y_{i},\epsilon_{i})$ are pairwise disjoint subsets of $N$, and
choose $h_{i}\in\Cinf(N)$ with compact support in
$K(y_{i},\epsilon_{i})$ such that
\[ (h_{i})_{y_{i}}=1 \]
for any $i$. By the definition of $R$ we can find for every $i$
a polynomial $P_{i}\in \Pi_{r}$
such that
\[ \left| D^{\alpha}(h_{i}P_{i})(y) \right| < 1/i \]
for any $y\in N$ and any multi-index $\alpha$ of order
$|\alpha|\leq r$, and
\[ \calA (P_{i})(x_{i}) \geq i \;.\]
These assumptions imply that the sum
\[ h=\sum_{i} h_{i}P_{i} \]
and all its partial derivatives of order at most $r$
converge uniformly on $N$,
thus $h\in\Cr(N)$. On the other hand, by Corollary \ref{corollary4}
we have
\[ \calA (h)(x_{i})=\calA (P_{i})(x_{i}) \geq i \;, \]
which contradicts the continuity of $\calA(h)$ at $x$.
The set $M\setminus R$ is therefore closed and discrete in $M$.

Take any $x\in R$. By \cite[page 35]{AczelDhombres}
it follows that the additive map $\Phi_{x}$ is in fact
linear. If we linearly parameterise $\Pi_{r}$ by
the partial derivatives of the polynomials in $\Pi_{r}$
at $\phi(x)$,
we obtain unique real numbers
$p_{\alpha}(x)$, for $|\alpha|\leq r$,
such that
\begin{equation} \label{equation31}
\Phi_{x}(P) = \calA(P)(x)
   = \sum_{|\alpha|\leq r} p_{\alpha}(x) D^{\alpha}(P)(\phi(x))
\end{equation}
for any $P\in \Pi_{r}$.
By induction on $|\alpha|$ we can
check that all the functions $p_{\alpha}$
are continuous on $R$. Indeed, 
if we take $P$ to be the homogeneous polynomial
$P(t_{1},\ldots,t_{n})=t_{1}^{\alpha_{1}}\cdots t_{n}^{\alpha_{n}}$,
$\alpha=(\alpha_{1},\ldots,\alpha_{n})$,
we obtain from (\ref{equation31})
an explicit polynomial expression for $p_{\alpha}(x)$
in terms of
$\calA(P)(x)$, $\phi(x)$ and $p_{\beta}(x)$, for $|\beta|<|\alpha|$.

Take any $g\in\Cr(N)$. Let $P_{g}\in\Pi_{r}$ be the
Taylor polynomial of $g$ of order $r$ around $\phi(x)$, i.e.\
$j^{r}_{\phi(x)}(P_{g})=j^{r}_{\phi(x)}(g)$. It follows from
Lemma \ref{lemma5} that
$\calA(P_{g})(x)=\calA(g)(x)$. The equation
(\ref{equation31}) therefore implies
\[ \calA(g)(x)
   = \sum_{|\alpha|\leq r} p_{\alpha}(x) D^{\alpha}(g)(\phi(x)) \]
for any $g\in\Cr(N)$ and any $x\in R$.
\end{proof}

\begin{lemma} \label{lemma7}
Suppose that $r<\infty$. Then
\[ \calA(g)(x) = \calA(1)(x) g(\phi(x)) \]
for any $g\in \Cr(N)$ and any $x\in M$.
\end{lemma}

\begin{proof}
By Corollary \ref{corollary4} we can assume without
loss of generality that $M$ and $N$ are open subsets
of $\RR^{n}$. By Lemma \ref{lemma6}
there exist an open subset $R$ with
discrete complement in $M$
and continuous real functions $p_{\alpha}$ on $R$,
for any multi-index $\alpha$ of order $|\alpha|\leq r$,
such that
\begin{equation} \label{equation33}
\calA(g)(x) = \sum_{|\alpha|\leq r} p_{\alpha}(x) D^{\alpha}(g)(\phi(x))
\end{equation}
for any $g\in \Cr(N)$ and any $x\in R$.
Observe that $p_{(0,\ldots,0)}=\calA(1)$.

Let $I$ be the set of all multi-indices $\alpha$ of order
$|\alpha|\leq r$ such that $p_{\alpha}$ is not identically
zero on $R$. Note that $I\neq\emptyset$ because $\calA\neq 0$.
Choose $\alpha\in I$ with a component
$\alpha_{i}$ which is maximal among all the components of all
multi-indices in $I$, i.e.\ $\alpha_{i}\geq \beta_{j}$ for any
$\beta\in I$ and any $1\leq j\leq n$. We will show that $\alpha_{i}=0$,
which implies $I=\{(0,\ldots,0)\}$.

Suppose that $\alpha_{i}>0$.
Observe that $k=r+1-\alpha_{i}$ satisfies $0<k\leq r$.
Since $p_{\alpha}$ is not identically zero on $R$, there exists
an open non-empty connected
subset $U$ of $R$ such that $p_{\alpha}$ has
no zeros on $U$. Take any $a\in U$ and choose $\epsilon>0$
so small that $K(a,\epsilon)\subset U$.
Define a multi-index $\alpha'$ by
$\alpha'_{i}=r+1$ and $\alpha'_{l}=\alpha_{l}$ for any $l\neq i$.

Let $P\in\Cr(N)$ be the homogeneous polynomial function given
by
\[ P(t) = (t-\phi(a))^{\alpha'} =
   (t_{1}-\phi_{1}(a))^{\alpha'_{1}}\cdots
   (t_{n}-\phi_{n}(a))^{\alpha'_{n}} \;,
\]
where $t=(t_{1},\ldots,t_{n})\in N\subset\RR^{n}$
and $\phi(x)=(\phi_{1}(x),\ldots,\phi_{n}(x))$
for any $x\in M$.

Take any $\beta\in I$. We write
$\beta\leq\alpha'$ if $\beta_{l}\leq\alpha'_{l}$
for all $l$. For the derivative
$D^{\beta}(P)$ we have the following two possibilities:

(i) If $\beta\not\leq\alpha'$, then $D^{\beta}(P)=0$.

(ii) If $\beta\leq\alpha'$, then we have in fact 
$\beta\leq\alpha$ (because $\alpha_{i}\geq\beta_{i}$ by the maximality
of the component $\alpha_{i}$ of $\alpha$)
and
\[ D^{\beta}(P)(t)= c_{\beta} \,(t-\phi(a))^{\alpha'-\beta}
   = c_{\beta} \,(t-\phi(a))^{\alpha-\beta} (t_{i}-\phi_{i}(a))^{k} \]
for some non-zero $c_{\beta}\in\RR$.
The polynomial
$c_{\beta} \,(t-\phi(a))^{\alpha-\beta}$
is constant and non-zero only in the case $\beta=\alpha$.
In all other cases it has value $0$ at $\phi(a)$.
Put $I'=\{ \beta\in I \,|\, \beta\leq\alpha\,,\; \beta\neq\alpha \}$
and denote
\[ w(t) = \sum_{\beta\in I'} c_{\beta}\,p_{\beta}(\phi^{-1}(t))
          (t-\phi(a))^{\alpha-\beta} \;. \]
This is a continuous function of $t\in \phi(R)$
which equals $0$ at $\phi(a)$.

Note that $j^{r}_{\phi(a)}(P)=j^{r}_{\phi(a)}(0)$,
so Lemma \ref{lemma5} implies that
$j^{r}_{a}(\calA(P))=j^{r}_{a}(0)$.
Choose any $j=1,\ldots,n$,
and denote by $e_{j}$ the $j$-th vector of the standard basis
of $\RR^{n}$.
The Taylor formula for $\calA(P)$ at $a$
gives for any $|h|<\epsilon$ a real number
$0<\vartheta<1$ such that
\begin{equation} \label{equation34}
\calA(P)(a+he_{j})=
   \frac{1}{k!}\frac{\partial^{k}\calA(P)}
        {\partial t_{j}^{k}}(a+\vartheta h e_{j})
   h^{k} \;.
\end{equation}
Note that the function
\[ z(h)=
\frac{1}{k!}\frac{\partial^{k}\calA(P)}
                 {\partial t_{j}^{k}}(a+\vartheta h e_{j}) \]
is continuous in $h=0$ and satisfies $z(0)=0$.

On the other hand, from (\ref{equation33}) it follows that
\begin{equation} \label{equation43}
\begin{split}
\calA(P)(a+he_{j})
   &= \sum_{|\beta|\leq r} p_{\beta}(a+he_{j}) D^{\beta}(P)(\phi(a+he_{j})) \\
   &= \sum_{\beta\in I} c_{\beta}p_{\beta}(a+he_{j})
       (\phi(a+he_{j})-\phi(a))^{\alpha'-\beta} \\
   &= \left( w(\phi(a+he_{j})) + c_{\alpha}\,p_{\alpha}(a+he_{j}) \right)
      (\phi_{i}(a+he_{j})-\phi_{i}(a))^{k}
\end{split}
\end{equation}
By combining (\ref{equation34}) and (\ref{equation43}) we obtain
for $0<|h|<\epsilon$
\[
\left( w(\phi(a+he_{j})) + c_{\alpha}\,p_{\alpha}(a+he_{j}) \right)
\left( \frac{\phi_{i}(a+he_{j})-\phi_{i}(a)}{h} \right)^{k} = z(h)\;.
\]
When $h$ approaches $0$, the right hand side converges to $0$,
but the first factor of the left hand side converges to
$c_{\alpha}\,p_{\alpha}(a)\neq 0$, which is possible only if the limit
\[
\lim_{h\ra 0}\frac{\phi_{i}(a+he_{j})-\phi_{i}(a)}{h}
= \frac{\partial\phi_{i}}{\partial t_{j}}(a)
\]
exists and equals $0$.
Since this is true for any $j$ and any point $a\in U$,
it follows that $\phi_{i}$ is constant on $U$. In particular,
the restriction $\phi|_{U}$ is not open, which is in
contradiction with
the fact that $\phi$ is a homeomorphism.

We can therefore conclude that $I=\{(0,\ldots,0)\}$ and hence
\[ \calA(g)(x) = \calA(1)(x) g(\phi(x)) \]
for any $g\in\Cr(N)$ and any $x\in R$. Because $R$
is dense in $M$ and both sides of the last equation are continuous
functions of $x$, defined on all of $M$, it follows that
the equality holds true for any $x\in M$.
\end{proof}

\begin{proof}[Proof of Theorem \ref{theorem1}]
For any positive $g\in\Cr(N)$ we have
$\calB(g)={\exp}\com \calA({\ln}\com g)$, so 
Lemma \ref{lemma7} implies that
\begin{equation} \label{equation11}
\calB(g)(x) = g(\phi(x))^{p(x)} 
\end{equation}
for any $x\in M$,
where $p=\calA(1)={\ln}\com\calB(\mathrm{e}1)\in\Cr(M)$.
Since the inverse $\phi^{-1}$ corresponds to the multiplicative
bijection $\calB^{-1}$, it follows analogously from
Lemma \ref{lemma7} that there exists $q\in\Cr(N)$
such that
\begin{equation} \label{equation13}
   \calB^{-1}(f)(y)=f(\phi^{-1}(y))^{q(y)}
\end{equation}
for any positive $f\in\Cr(M)$ and for
any $y\in N$.
Direct computation of the composition of $\calB$ and $\calB^{-1}$
shows that $(p\com\phi^{-1})q=1$, so in particular both $p$ and $q$
are nowhere zero. 
Furthermore, it follows from (\ref{equation11}) that
the composition of a positive $\Cr$-function with $\phi$
is again a $\Cr$-function, thus $\phi$ is of class $\Cr$.
Analogously, by (\ref{equation13}) the map $\phi^{-1}$ is
of class $\Cr$, thus we may conclude that $\phi$ is
a $\Cr$-diffeomorphism.

We will now show that both $p=1$ and $q=1$.
Because $(p\com\phi^{-1})q=1$, it is sufficient to show
that $p\geq 1$ and $q\geq 1$. By symmetry it is enough to
prove $p\geq 1$. 
So assume that there is a point $x\in M$
such that $p(x)<1$.
Since $\phi$ is a $\Cr$-diffeomorphism, we may choose a $\Cr$-path
$\gamma\!:(-1,1)\ra M$ with $\gamma(0)=x$
and a function $g\in\Cr(N)$
such that $g(\phi(\gamma(t)))=t$ for any $t\in (-1,1)$.
Note that $\sigma=p\com\gamma$
and $u=\calB(g)\com \gamma$ are $\Cr$-functions on $(-1,1)$.
It follows from (\ref{equation11}) and Lemma \ref{lemma2} (i) that
\[ u(t)=t^{\sigma(t)} \]
for any $t>0$. Since $u$ is continuous on $(-1,1)$,
this is possible only if $\sigma(0)\geq 0$.
As $p$ has no zeros, this implies $p(x)=\sigma(0)>0$. 
By the continuity of $\sigma$ we may choose
$0<a<b<1$ and $0<\epsilon<1$ such that
$a<\sigma(t)<b$ for any $t\in(-\epsilon,\epsilon)$.

The derivative of $u$ at a point $t\in(0,1)$ equals
\begin{equation} \label{equation12}
\frac{t^{\sigma(t)} \sigma(t) }{t}
+ t^{\sigma(t)} \ln(t) \frac{d\sigma}{dt}(t)\;.
\end{equation}
Since the derivative of $\sigma$
is bounded on a small neighbourhood of $0$ and
$\sigma(-\epsilon,\epsilon)\subset [a,b]\subset (0,1)$,
it follows that the second summand
of (\ref{equation12})
converges to $0$ as $t>0$ approaches $0$.
On the other hand, the fact that
$\sigma(-\epsilon,\epsilon)\subset [a,b]\subset (0,1)$
implies that the first summand of (\ref{equation12})
is unbounded on any neighbourhood
of $0$. Hence $u\in\Cr(-1,1)$ has unbounded derivative
on any neighbourhood of $0$,
which is a contradiction.

Thus $p=1$, and therefore
\[ \calB(g)(x)=g(\phi(x)) \]
for any positive $g\in\Cr(N)$ and
any $x\in M$. Finally, it follows from
Lemma \ref{lemma2} that this formula actually holds
for any $g\in \Cr(N)$.
\end{proof}

\end{document}